\magnification=1200

\hsize=11.25cm    
\vsize=18cm       
\parindent=12pt   \parskip=5pt     

\hoffset=.5cm   
\voffset=.8cm   

\pretolerance=500 \tolerance=1000  \brokenpenalty=5000

\catcode`\@=11

\font\eightrm=cmr8         \font\eighti=cmmi8
\font\eightsy=cmsy8        \font\eightbf=cmbx8
\font\eighttt=cmtt8        \font\eightit=cmti8
\font\eightsl=cmsl8        \font\sixrm=cmr6
\font\sixi=cmmi6           \font\sixsy=cmsy6
\font\sixbf=cmbx6

\font\tengoth=eufm10 
\font\eightgoth=eufm8  
\font\sevengoth=eufm7      
\font\sixgoth=eufm6        \font\fivegoth=eufm5

\skewchar\eighti='177 \skewchar\sixi='177
\skewchar\eightsy='60 \skewchar\sixsy='60

\newfam\gothfam           \newfam\bboardfam

\def\tenpoint{
  \textfont0=\tenrm \scriptfont0=\sevenrm \scriptscriptfont0=\fiverm
  \def\rm{\fam\z@\tenrm}
  \textfont1=\teni  \scriptfont1=\seveni  \scriptscriptfont1=\fivei
  \def\oldstyle{\fam\@ne\teni}\let\old=\oldstyle
  \textfont2=\tensy \scriptfont2=\sevensy \scriptscriptfont2=\fivesy
  \textfont\gothfam=\tengoth \scriptfont\gothfam=\sevengoth
  \scriptscriptfont\gothfam=\fivegoth
  \def\goth{\fam\gothfam\tengoth}
  
  \textfont\itfam=\tenit
  \def\it{\fam\itfam\tenit}
  \textfont\slfam=\tensl
  \def\sl{\fam\slfam\tensl}
  \textfont\bffam=\tenbf \scriptfont\bffam=\sevenbf
  \scriptscriptfont\bffam=\fivebf
  \def\bf{\fam\bffam\tenbf}
  \textfont\ttfam=\tentt
  \def\tt{\fam\ttfam\tentt}
  \abovedisplayskip=12pt plus 3pt minus 9pt
  \belowdisplayskip=\abovedisplayskip
  \abovedisplayshortskip=0pt plus 3pt
  \belowdisplayshortskip=4pt plus 3pt 
  \smallskipamount=3pt plus 1pt minus 1pt
  \medskipamount=6pt plus 2pt minus 2pt
  \bigskipamount=12pt plus 4pt minus 4pt
  \normalbaselineskip=12pt
  \setbox\strutbox=\hbox{\vrule height8.5pt depth3.5pt width0pt}
  \let\bigf@nt=\tenrm       \let\smallf@nt=\sevenrm
  \normalbaselines\rm}

\def\eightpoint{
  \textfont0=\eightrm \scriptfont0=\sixrm \scriptscriptfont0=\fiverm
  \def\rm{\fam\z@\eightrm}
  \textfont1=\eighti  \scriptfont1=\sixi  \scriptscriptfont1=\fivei
  \def\oldstyle{\fam\@ne\eighti}\let\old=\oldstyle
  \textfont2=\eightsy \scriptfont2=\sixsy \scriptscriptfont2=\fivesy
  \textfont\gothfam=\eightgoth \scriptfont\gothfam=\sixgoth
  \scriptscriptfont\gothfam=\fivegoth
  \def\goth{\fam\gothfam\eightgoth}
  
  \textfont\itfam=\eightit
  \def\it{\fam\itfam\eightit}
  \textfont\slfam=\eightsl
  \def\sl{\fam\slfam\eightsl}
  \textfont\bffam=\eightbf \scriptfont\bffam=\sixbf
  \scriptscriptfont\bffam=\fivebf
  \def\bf{\fam\bffam\eightbf}
  \textfont\ttfam=\eighttt
  \def\tt{\fam\ttfam\eighttt}
  \abovedisplayskip=9pt plus 3pt minus 9pt
  \belowdisplayskip=\abovedisplayskip
  \abovedisplayshortskip=0pt plus 3pt
  \belowdisplayshortskip=3pt plus 3pt 
  \smallskipamount=2pt plus 1pt minus 1pt
  \medskipamount=4pt plus 2pt minus 1pt
  \bigskipamount=9pt plus 3pt minus 3pt
  \normalbaselineskip=9pt
  \setbox\strutbox=\hbox{\vrule height7pt depth2pt width0pt}
  \let\bigf@nt=\eightrm     \let\smallf@nt=\sixrm
  \normalbaselines\rm}

\tenpoint

\def\pc#1{\bigf@nt#1\smallf@nt}         \def\pd#1 {{\pc#1} }

\catcode`\;=\active
\def;{\relax\ifhmode\ifdim\lastskip>\z@\unskip\fi
\kern\fontdimen2  -1.2 \fontdimen3 \string;}

\catcode`\:=\active
\def:{\relax\ifhmode\ifdim\lastskip>\z@\unskip\fi\penalty\@M\ \fi\string:}

\catcode`\!=\active
\def!{\relax\ifhmode\ifdim\lastskip>\z@
\unskip\fi\kern\fontdimen2  -1.1 \fontdimen3 \string!}

\catcode`\?=\active
\def?{\relax\ifhmode\ifdim\lastskip>\z@
\unskip\fi\kern\fontdimen2  -1.1 \fontdimen3 \string?}

% \catcode`\«=\active 
% \def«{\raise.4ex\hbox{%
%  $\scriptscriptstyle\langle\!\langle$}}
% 
% \catcode`\»=\active 
% \def»{\raise.4ex\hbox{%
%  $\scriptscriptstyle\rangle\!\rangle$}}

\frenchspacing

\def\raggedbottom{\topskip 10pt plus 36pt\r@ggedbottomtrue}

\def\pointir{\unskip . --- \ignorespaces}

\def\Medbreak{\vskip-\lastskip\medbreak}

\long\def\th#1 #2\enonce#3\endth{
   \Medbreak\noindent
   {\pc#1} {#2\unskip}\pointir{\it #3}\smallskip}

\def\proof{\vskip-\lastskip\smallskip\noindent
 {\it Proof} : }

\def\decale#1{\smallbreak\hskip 28pt\llap{#1}\kern 5pt}
\def\decaledecale#1{\smallbreak\hskip 34pt\llap{#1}\kern 5pt}
\def\puce{\smallbreak\hskip 6pt{$\scriptstyle\bullet$}\kern 5pt}

\def\eqalign#1{\null\,\vcenter{\openup\jot\m@th\ialign{
\strut\hfil$\displaystyle{##}$&$\displaystyle{{}##}$\hfil
&&\quad\strut\hfil$\displaystyle{##}$&$\displaystyle{{}##}$\hfil
\crcr#1\crcr}}\,}

\catcode`\@=12

\showboxbreadth=-1  \showboxdepth=-1

\newcount\numerodesection \numerodesection=1
\def\section#1{\bigbreak
 {\bf\number\numerodesection.\ \ #1}\nobreak\medskip
 \advance\numerodesection by1}

\mathcode`A="7041 \mathcode`B="7042 \mathcode`C="7043 \mathcode`D="7044
\mathcode`E="7045 \mathcode`F="7046 \mathcode`G="7047 \mathcode`H="7048
\mathcode`I="7049 \mathcode`J="704A \mathcode`K="704B \mathcode`L="704C
\mathcode`M="704D \mathcode`N="704E \mathcode`O="704F \mathcode`P="7050
\mathcode`Q="7051 \mathcode`R="7052 \mathcode`S="7053 \mathcode`T="7054
\mathcode`U="7055 \mathcode`V="7056 \mathcode`W="7057 \mathcode`X="7058
\mathcode`Y="7059 \mathcode`Z="705A

% handling accented characters in plain TeX :

% \catcode`\À=\active \defÀ{\`A}    \catcode`\à=\active \defà{\`a} 
% \catcode`\Â=\active \defÂ{\^A}    \catcode`\â=\active \defâ{\^a} 
% \catcode`\Æ=\active \defÆ{\AE}    \catcode`\æ=\active \defæ{\ae}
% \catcode`\Ç=\active \defÇ{\c C}   \catcode`\ç=\active \defç{\c c}
% \catcode`\È=\active \defÈ{\`E}    \catcode`\è=\active \defè{\`e} 
% \catcode`\É=\active \defÉ{\'E}    \catcode`\é=\active \defé{\'e} 
% \catcode`\Ê=\active \defÊ{\^E}    \catcode`\ê=\active \defê{\^e} 
% \catcode`\Ë=\active \defË{\"E}    \catcode`\ë=\active \defë{\"e} 
% \catcode`\Î=\active \defÎ{\^I}    \catcode`\î=\active \defî{\^\i}
% \catcode`\Ï=\active \defÏ{\"I}    \catcode`\ï=\active \defï{\"\i}
% \catcode`\Ô=\active \defÔ{\^O}    \catcode`\ô=\active \defô{\^o} 
% \catcode`\Ù=\active \defÙ{\`U}    \catcode`\ù=\active \defù{\`u} 
% \catcode`\Û=\active \defÛ{\^U}    \catcode`\û=\active \defû{\^u} 
% \catcode`\Ü=\active \defÜ{\"U}    \catcode`\ü=\active \defü{\"u} 

\def\mod{\mathop{\rm mod.}\nolimits}
\def\pmod#1{\;(\mod#1)}

\openup1\jot

\font\tencyr=wncyr10
\font\sevencyr=wncyr7
\font\fivecyr=wncyr5
\newfam\cyrfam
\def\cyr{\fam\cyrfam\tencyr}
\textfont\cyrfam=\tencyr
\scriptfont\cyrfam=\sevencyr
\scriptscriptfont\cyrfam=\fivecyr

\newcount\refno

\long\def\ref#1:#2<#3>{                                        
\global\advance\refno by1\par\noindent                              
\llap{[{\bf\number\refno}]\ }{#1}\pointir{\it #2} #3\goodbreak }

\def\citer#1(#2){[{\bf\number#1}\if#2\empty\relax\else,\ #2\fi]}

\newbox\bibbox
\setbox\bibbox\vbox{
\bigskip
\centerline{---$*$---$*$---}
\bigbreak
\centerline{{\pc BIBLIOGRAPHICAL REFERENCES}}

\def\mr#1.{{\sevenrm MR#1.}}

\ref{\pc CASSELS} (John):
Diophantine equations with special reference to elliptic curves,
<J.\ London Math.\ Soc.\ {\bf 41}, 1966, pp.~193--291.>
\newcount\cassels  \global\cassels=\refno

\ref {\pc COATES} (John) and {\pc WILES} (Andrew) :
On the conjecture of Birch and Swinnerton-Dyer,
<Invent.\ Math.\ {\bf 39} (1977), no.\ 3, 223--251.> % \mr0463176
\newcount\coateswiles  \global\coateswiles=\refno

\ref {\pc ELKIES} (Noam) :
The existence of infinitely many supersingular primes for every elliptic curve
over $Q$,
<Invent.\ Math.\ {\bf 89} (1987), no.\ 3, 561--567.> % \mr0903384
\newcount\elkies  \global\elkies=\refno

\ref {\pc FALTINGS} (Gerd) :
Endlichkeitss{\"a}tze f{\"u}r abelsche Variet{\"a}ten {\"u}ber
Zahlk{\"o}rpern,
<Invent.\ Math.\ {\bf 73} (1983), no.\ 3, 349--366.>  % \mr0718935
\newcount\faltings \global\faltings=\refno

\ref {\pc GROSS} (Benedict) and {\pc ZAGIER} (Don) :
Heegner points and derivatives of\/ $L$-series,
<Invent.\ Math.\ {\bf 84} (1986), no.\ 2, 225--320.> % \mr0833192
\newcount\grosszagier \global\grosszagier=\refno

\ref{\pc HEEGNER} (Kurt) :
Diophantische Analysis und Modulfunktionen,
<Math.\ Z.\ {\bf 56} (1952), 227--253.> % \mr0053135
\newcount\heegner \global\heegner=\refno

\ref{\pc KATO} (Kazuya) and {\pc TRIHAN} (Fabien):
On the conjectures of Birch and Swinnerton-Dyer in characteristic $p>0$,
< Invent.\ Math.\  {\bf 153}  (2003),  no.\ 3, 537--592.> % \mr2000469
\newcount\katotrihan \global\katotrihan=\refno

\ref {\pc KOLYVAGIN} (Victor) :
The Mordell-Weil and Shafarevich-Tate groups for Weil elliptic curves,
<(Russian) Izv.\ Akad.\ Nauk SSSR Ser.\ Mat.\ {\bf 52} (1988), no.\ 6,
1154--1180, 1327.> % \mr0984214
\newcount\kolyvagin  \global\kolyvagin=\refno

\ref {\pc LANGLANDS} (Robert) :
Problems in the theory of automorphic forms,
<Lecture Notes in Math., Vol.\ 170, Springer, Berlin, 1970,
pp. 18--61.>  % \mr0302614
\newcount\langlands  \global\langlands=\refno

\ref {\pc MAZUR} (Barry) :
Modular curves and the Eisenstein ideal,
<Inst.\ Hautes {\'E}tudes Sci.\ Publ.\ Math.\ {\bf 47} (1977), 33--186
(1978).> % \mr0488287
\newcount\eisenstein  \global\eisenstein=\refno

\ref {\pc MAZUR} (Barry) :
On the passage from local to global in number theory,
<Bull.\ Amer.\ Math.\ Soc.\ (N.S.) {\bf 29} (1993), no. 1, 14--50.> 
% \mr1202293   
\newcount\mazur  \global\mazur=\refno

\ref {\pc  RUBIN} (Karl) :
The ``main conjectures'' of Iwasawa theory for imaginary quadratic fields,
<Invent.\ Math.\ {\bf 103} (1991), no.\ 1, 25--68.> % \mr1079839
\newcount\rubin  \global\rubin=\refno

\ref {\pc SERRE} (Jean-Pierre) :
Propri{\'e}t{\'e}s galoisiennes des points d'ordre fini des courbes
elliptiques, 
<Invent.\ Math.\ {\bf 15} (1972), no.\ 4, 259--331.> % \mr0387283
\newcount\serre  \global\serre=\refno

\ref {\pc SERRE} (Jean-Pierre) and {\pc TATE} (John) :
Good reduction of abelian varieties,
<Ann.\ of Math.\ (2) {\bf 88} (1968) 492--517.> % \mr0236190
\newcount\serretate  \global\serretate=\refno

\ref {\pc SKINNER} (Christopher) and {\pc URBAN} ({\'E}ric) :
Sur les d{\'e}formations $p$-adiques de certaines repr{\'e}sentations
automorphes,
<J.\ Inst.\ Math.\ Jussieu {\bf 5} (2006), no.\ 4, 629--698.> % \mr2261226
\newcount\skinnerurban \global\skinnerurban=\refno

\vfill\eject

\ref {\pc TATE} (John) :
The arithmetic of elliptic curves,
<Invent.\ Math.\ {\bf 23} (1974), 179--206.> % \mr0419359
\newcount\tate \global\tate=\refno

\ref {\pc TAYLOR} (Richard) :
Automorphy for some $l$-adic lifts of automorphic mod~$l$ representations. II, 
<{\tt www.math.harvard.edu/\char`\~rtaylor/}.>
\newcount\taylor \global\taylor=\refno

\ref {\pc TUNNELL} (Jerrold) :
A classical Diophantine problem and modular forms of weight $3/2$,
<Invent.\ Math.\ {\bf 72} (1983), no.\ 2, 323--334.> % \mr0700775
\newcount\tunnell  \global\tunnell=\refno

\medskip 

The site {\tt gdz.sub.uni-goettingen.de} has the papers \citer\coateswiles(),
\citer\elkies(), \citer\faltings(), \citer\grosszagier(), \citer\heegner(),
\citer\rubin(), \citer\serre(), \citer\tate() and \citer\tunnell().  The site
{\tt www.numdam.org} has the paper \citer\eisenstein(). 
All papers by Langlands are available at his website.
 
}

\centerline{\bf Congruent numbers, elliptic curves,}
\centerline{\bf and the passage from the local to the
global} 

\vskip1cm
\centerline{Chandan Singh Dalawat}
\medskip

{\it The ancient unsolved problem of congruent numbers has been reduced to one
  of the major questions of contemporary arithmetic~: the finiteness of the
  number of curves over $\bf Q$ which become isomorphic at every place to a
  given curve.  We give an elementary introduction to congruent numbers and
  their conjectural characterisation, discuss local-to-global issues leading
  to the finiteness problem, and list a few results and conjectures in the
  arithmetic theory of elliptic curves.}

\bigskip

The area $\alpha$ of a right triangle with sides $a,b,c$ (so that
$a^2+b^2=c^2$) is given by $2\alpha=ab$.  If $a,b,c$ are rational, then so is
$\alpha$.  Conversely, which rational numbers $\alpha$ arise as the area of a
rational right triangle $a,b,c$~?  This problem of characterising ``congruent
numbers'' --- areas of rational right triangles --- is perhaps the oldest
unsolved problem in all of Mathematics.  It dates back to more than a thousand
years and has been variously attributed to the Arabs, the Chinese, and the
Indians.

Three excellent accounts of the problem are available on the Web~: {\it Right
  triangles and elliptic curves\/} by Karl Rubin, {\it Le probl{\`e}me des
  nombres congruents\/} by Pierre Colmez, which also appears in the
October~2006 issue of the {\it Gazette des math{\'e}maticiens}, and Franz
Lemmermeyer's translation {\it Congruent numbers, elliptic curves, and modular
  forms} of an article in French by Guy Henniart.  A more elementary
introduction is provided by the notes of a lecture in Hong Kong by John
Coates, which have appeared in the August~2005 issue of the {\it Quaterly
  journal of pure and applied mathematics}.  A detailed account is to be found
in the {\it Introduction to elliptic curves and modular forms\/} (Springer,
1984) by Neal Koblitz.  None of these sources goes beyond the theorems of
Coates \& Andrew Wiles \citer\coateswiles() (see Theorem~14) and of Jerrold
Tunnell \citer\tunnell() (see Theorem~25).
  
  In 1991, Rubin \citer\rubin() (see Theorem~15) reduced the congruent number
  problem to a natural finiteness question in the arithmetic of elliptic
  curves (with ``complex multiplications'').  An excellent survey of such
  finiteness questions can be found in Barry Mazur's article \citer\mazur().
  
  These notes consist of three parts of quite different nature.  The first
  part is an elementary presentation of the problem of congruent numbers (\S1)
  and its conjectural solution (\S2)~; the material here is borrowed from the
  accounts which have been cited.  The second part introduces local number
  fields (\S3) and discusses the local-to-global principle --- its validity in
  the case of conics (\S4) and its failure in the case of cubics (\S5) --- in
  a language which can be understood by bright undergraduates.  The last part,
  which requires greater mathematical maturity, is a catalogue of results ---
  some old, some new --- and conjectures in the arithmetic theory of elliptic
  curves in general (\S6) and those without complex multiplications in
  particular (\S7)~; it ends with a word about the role played by modular
  forms (\S8).

{\eightpoint I thank Pere Clark for a very careful reading of the manuscript,
  and for his suggestions for improvement.}

\bigbreak 
\centerline{\bf 1. Congruent numbers} 
\medskip 

If a rational number $\alpha$ is the area of a right triangle with rational
sides, then so is $\alpha\beta^2$ for every rational $\beta\in{\bf Q}^\times$.
Indeed, if $\alpha$ is the area of a rational right triangle with sides
$a,b,c$, then $\alpha\beta^2$ is the area of the rational right triangle with
sides $a|\beta|,b|\beta|,c|\beta|$.  So, up to replacing $\alpha$ by
$\alpha\beta^2$ for a suitable $\beta$, we may assume that $\alpha$ is an
integer, and moreover that $\alpha$ is not divisible by the square of any
prime number.  In other words, we assume that $\alpha$ is a positive squarefree
integer.

\th DEFINITION 1
\enonce
A positive squarefree integer $\alpha$ is said to be a congruent number if
there exist $a,b,c\in{\bf Q}$ such that $a^2+b^2=c^2$ and $ab=2\alpha$.
\endth

The terminology is classical and comes from the fact that $\alpha$ is
congruent if and only if it is the common difference ({\it congruum}, in
Latin) of a three-term arithmetic progression of rational squares.  For if
$\alpha$ is the area of a rational right triangle with sides $a<b<c$, then,
putting $d=(c/2)^2$, the arithmetic progression $d-\alpha$, $d$, $d+\alpha$
consists of rational squares.  Conversely, if there is a rational number $d$
such that $d-\alpha$, $d$, $d+\alpha$ are all three squares, then $\alpha$ is
the area of the rational right triangle with sides
$\sqrt{d+\alpha}-\sqrt{d-\alpha}$, $\sqrt{d+\alpha}+\sqrt{d-\alpha}$ and
$2\sqrt d$.  

The problem is to determine which numbers are congruent.  Let us first study
the single equation $a^2+b^2=c^2$ in rational numbers $>0$.  Two solutions
$(a,b,c)$, $(a',b',c')$ of this equation are called {\it equivalent\/} if
$a=\lambda.a'$, $b=\lambda.b'$, $c=\lambda.c'$ for some $\lambda\in{\bf
  Q}^\times$.  A rational solution is called {\it primitive\/} if
$a,b,c\in{\bf Z}$, and if they have no common prime divisor.  Every rational
solution is equivalent to a primitive one, and no two primitive solutions are
equivalent. 

Reducing a primitive solution modulo~$4$, we see that precisely one of $a,b$
is even.

\th PROPOSITION 2
\enonce
Let $(a,b,c)$ be a primitive solution of $a^2+b^2=c^2$, with $a=2t$ even.
Then there exist integers $m>n>0$, $\gcd(m,n)=1$, $m\not\equiv n\pmod2$, such
that 
$$
a=2mn,\ \ b=m^2-n^2,\ \ c=m^2+n^2.\leqno{(1)}
$$
\endth \proof As $b$ is odd, so is $c$.  Hence $c+b$ and $c-b$ are even;
write $c+b=2u$ and $c-b=2v$.  If a prime number divides both $u$ and $v$, it
would divide their sum $u+v=c$ and their difference $u-v=b$.  But
$\gcd(b,c)=1$, so we have $\gcd(u,v)=1$.  The relation $a^2+b^2=c^2$ implies
that $t^2=uv$, which shows that each of $u,v$ must be a square.  Putting
$u=m^2$, $v=n^2$ proves $(1)$.  Finally, if $m\equiv n\pmod2$, then $2$ would
divide $a,b,c$ and the solution would not be primitive.

% {\eightpoint This proposition can also be proved by determining the rational
%   points on the circle $C:x^2+y^2=1$ obtained upon dividing throughout by
%   $c^2$.  The point $O=(1,0)$ lies on $C$; let us fix a rational line $L$ not
%   passing through $O$.  Rational lines through $O$ meet the circle $C$ at a
%   rational point $P$ and the line $L$ at a rational point $f(P)$.  The map
%   $f$ gives a bijection between the sets $C({\bf Q})$, $L({\bf Q})$ of
%   rational points, if we include the points at infinity.}

{\eightpoint Let $C$ be a projective conic with a rational point $O$, for
  example the one defined by $a^2+b^2=c^2$, with $O=(1:0:1)$.  Denoting by $D$
  the projective line of lines through $O$, the morphism $f$ which sends a
  point $P\in C$ to the line $f(P)\in D$ passing through $O$ and $P$ --- the
  tangent to $C$ at $O$ if $P=O$ --- is an isomorphism.}

This result allows us to generate a list which will eventually contain any
given congruent number~: it suffices to go through the list of all such pairs
$(m,n)$, compute the area $mn(m^2-n^2)$ of the triangle $(1)$, and take the
``squarefree part''.  Thus the pair $(2,1)$ shows that the number
$6=2.1.(2^2-1^2)$ is congruent.

Retaining only the squarefree parts of the numbers produced by this procedure,
the first few congruent numbers which show up are
$$5,6,7,13,14,15,21,22,23,29,30,31,34,37,38,39,41,\ldots\leqno{(2)}$$

Note that we have {\it not\/} proved that the numbers $1,2,3$ are not
congruent; it may simply be that they haven't yet shown up on the list~!
Indeed, Leonardo of Pisa (called Fibonacci) (1175--1240) was challanged to
find a rational right triangle of area $5$ (he succeeded) and he conjectured
that $1$ is not congruent~; this was settled much later by Pierre Fermat
(1601--1665).

How can we determine if a specific number such as $157$ is congruent~?  The
na{\"\i ve} approach, suggested by the discussion just after Definition~1,
would be to go through a ``list'' of squares $d$ of rational numbers and to
see if both $d-157$ and $d+157$ are squares.  There is indeed such a
``list''~: first we go through the squares of the finitely many rational
numbers whose numerator and denominator have at most one digit, then through
the squares of those --- again finitely many --- whose numerator and
denominator have at most two digits, and so on.  It turns out that the first
square which works, according to Don Zagier, is
$$
d =\left({224403517704336969924557513090674863160948472041\over
2\times8912332268928859588025535178967163570016480830}\right)^{\!2}
\raise2pt\hbox{.}
$$
Clearly, this number could not have been found by the na{\"\i}ve approach~;
some theory is needed.  Also, as before, this approach cannot prove that the
given number, for example~$1$, is not congruent.

% The reader may wish to verify that the number $157$ does not show up on the
% list for a long time.  It is congruent nevertheless~; the ``smallest'' --- in
% the sense of having the smallest possible numerator and denominator ---
% rational right triangle with area $157$ turns out to be (Don Zagier)~:
% $$\eqalign{
% &a ={6803298487826435051217540\over 411340519227716149383203}
% \raise2pt\hbox{,}\qquad
% b ={411340519227716149383203\over 21666555693714761309610}
% \raise2pt\hbox{,}\cr
% &c ={224403517704336969924557513090674863160948472041\over
% 8912332268928859588025535178967163570016480830}
% \raise2pt\hbox{.}\cr}
% $$

\th THEOREM 3 (P. Fermat, $\sim1640$)
\enonce
The number $1$ is not congruent.
\endth
\proof 
We have to show that there is no rational right triangle whose area is a
square.  If there is such a triangle, we may assume, as before, that its sides
are integers not divisible by the square of any prime number.  Fermat's idea
of infinite descent consists in showing that if there were such a
``primitive'' triangle whose area is a square, then there would be a smaller 
primitive triangle whose area is also a square.  Clearly, this cannot go on
for ever.

Let $(a,b,c)$ be a primitive triangle whose area is a square.  Assume that $a$
is even and write $a=2mn$, $b=m^2-n^2$, $c=m^2+n^2$, with $\gcd(m,n)=1$
(Proposition~2).  As the area
$$
mn(m+n)(m-n)
$$ 
is a square, and as no two of the four factors have a common prime
divisor, all four must be squares:
$$m=x^2,\ \ n=y^2,\ \ m+n=u^2,\ \ m-n=v^2.$$
We have $\gcd(u,v)=1$, and both $u,v$ are odd because their product $b$ is
odd. We also have $u^2=v^2+2y^2$, which we rewrite as
$$
2y^2=(u+v)(u-v).\leqno{(3)}
$$
As $u,v$ are odd and $\gcd(u,v)=1$, we have $\gcd(u+v,u-v)=2$.  So one of
the two factors on the right in $(3)$ must be of the form $2r^2$ and the other
of the form $4s^2$.  In any case, the sum of their squares is $16s^4+4r^4$.
At the same time, $(u+v)^2+(u-v)^2=2(u^2+v^2)=4m=4x^2$.  Comparing these two
results, we get $4s^4+r^4=x^2$, which means that $(2s^2, r^2, x)$ is also an
integral right triangle whose area $(rs)^2$ is a square.  This triangle is
smaller than our original triangle $(a,b,c)$ because $x^4=m^2<m^2+n^2=c$~; it
may not be primitive, but the corresponding primitive triangle is even
smaller.

{\eightpoint The passage from the triple $(a,b,c)$ to the triple $(2s^2, r^2,
  x)$ can be construed as {\it division\/} by $\pm2$ on the elliptic curve
  $C_1:y^2=x^3-x$~; cf.~the discussion before Exercise~7, and the beginning of
  \S~6.  The idea of the size of a triple leads the notion of {\it height\/}
  of a rational point on an elliptic curves.}

\th COROLLARY 4
\enonce
The equation $x^4-y^4=z^2$ has no solutions in integers with  $xyz\neq0$.
\endth

\proof If there were a solution, the integral triangle $(2x^2y^2,
x^4-y^4, x^4+y^4)$ would have square area $(xyz)^2$.

\th COROLLARY 5
\enonce
The equation $x^4+y^4=z^4$ has no solutions in integers with $xyz\neq0$.
\endth
\proof If there were a solution, we would have $z^4-y^4=(x^2)^2$.

The system of equations ($a^2+b^2=c^2$; $ab=2\alpha$) whose solvability in
rational numbers characterises $\alpha$ as a congruent numbers can be changed
into a single, more familiar, equation.

\th PROPOSITION 6
\enonce  
The integer\/ $\alpha$ is congruent if and only if the equation  
$$
C_\alpha:\alpha y^2=x^3-x \leqno{(4)}
$$
has a solution $x,y\in{\bf Q}$ with $y\neq0$.  
\endth 
\proof If $(x,y)$ is such a solution, then the area of the rational right
triangle $(2|x|, |x^2-1|, |x^2+1|)$ is $\alpha$, up to a rational square
($y^2$).  Conversely, let $(a,b,c)$ be a rational right triangle and write
$$
a=\lambda.2mn,\ \ b=\lambda.(m^2-n^2),\ \ c=\lambda.(m^2+n^2)
\quad(m,n\in{\bf Z})
$$
for some $\lambda\in{\bf Q}^\times$ (Proposition~2).  If the area of this
triangle is $\alpha$, we have $\alpha=\lambda^2.mn(m^2-n^2)$, which means that
$(4)$ has the solution $x=m/n$ and $y=1/\lambda n$.

From a given rational point $P=(x,y)$ ($y\neq0$) on $C_\alpha$ $(4)$ we can
generate infinitely many others : the tangent to $C_\alpha$ at the point $P$
meets $C_\alpha$ at another rational point $P_1=(x_1,y_1)$, and this process
can be continued~; it can be shown not to terminate (cf.~the discussion of
torsion points on $C_\alpha$, before Theorem~18) .

\th EXERCISE 7
\enonce
$\displaystyle
(x_1,y_1)=
\left({(x^2+1)^2\over 4x(x^2-1)},-{x^6-5x^4-5x^2+1\over 8\alpha^3y^3}\right)$.
\endth

This has an amusing consequence which is not at all obvious at the outset~:

\th COROLLARY 8
\enonce
If (a squarefree positive integer) $\alpha$ is congruent, then it is the area
of infinitely many rational right triangles.
\endth

\th COROLLARY 9
\enonce
A squarefree positive integer $\alpha$ is congruent if and only if the
equation $\alpha y^2=x^3-x$ has infinitely many solutions $x,y\in{\bf Q}$. 
\endth

Here are the first few rational squares $d$ such that both $d-6$ and $d+6$ are
squares~: 
$$\def\\{\raise2pt\hbox{,}}
\left({5\over 2\times1}\right)^{\!2}\\
\left({1201\over 2\times70}\right)^{\!2}\\
\left({7776485\over 2\times1319901}\right)^{\!2}\\
\left({2094350404801\over 2\times241717895860}\right)^{\!2}\raise2pt\hbox{.}
$$

% Here are the three simplest rational right triangles of area $6$ :
% $$\def\\{\raise2pt\hbox{,}}
% (3,4,5),\ \left({7\over10}\\{120\over7}\\{1201\over70}\right),\ 
% \left({4653\over851}\\{3404\over1551}\\{7776485\over1319901}\right).
% $$

{\eightpoint The morphism $y\mapsto\sqrt\alpha y$ shows that the curves $C_1$
  and $C_\alpha$ become isomorphic over ${\bf Q}(\sqrt{\alpha})$~; this is
  expressed by saying that $C_\alpha$ is a ``quadratic twist'' of $C_1$.  The
  problem of congruent numbers thus consists in characterising the quadratic
  twists of the fixed elliptic curve $C_1$ which have infinitely many rational
  points.}

\bigbreak
\centerline{\bf 2. The conjectural solution}
\medskip

After these elementary observations, let us give the conjectural answer to the
problem of characterising congruent numbers.  

Recursively define the polynomial $g_r(T)=g_{r-1}(T)(1-T^{8r})(1-T^{16r})$,
starting with $g_1(T)=T(1-T^{8})(1-T^{16})$.  Notice that $g_r(T)-g_{r-1}(T)$
is of degree $>8r$, which means that the polynomials $g_r$ and $g_{r-1}$ have
the same terms till degree $8r$.  This implies that as $r\rightarrow+\infty$,
the $g_r$ tend to a formal series $g\in{\bf Z}[[T]]$.

\th NOTATION 10
\enonce
For\/ $j=1,2$ and integer\/ $n>0$, define\/ $c_j(n)$ as being the coefficient
of\/ $T^n$ in the formal series \/ $g(T)\theta_j(T)$, where  
$$
g(T)=T\prod_{n=1}^{+\infty}(1-T^{8n})(1-T^{16n})\quad\hbox{and}\quad
\theta_j(T)=1+2\sum_{n=1}^{+\infty}T^{2jn^2}. 
$$
\endth 

Notice that the numbers $c_j(n)$ are quite easy to compute.  Here are the
first few, for $n$ odd and squarefree.
$$\vbox{\halign{&\quad\hfil$#$\cr
\noalign{\vskip5pt\hrule\vskip4pt}
n&1&3&5&7&11&13&15&17&19&21&23\cr
\noalign{\vskip5pt\hrule\vskip4pt}
c_1(n)&1&2&0&0&-2&0&0&-4&-2&0&0\cr
\noalign{\vskip5pt\hrule\vskip4pt}
c_2(n)&1&0&2&0&0&-2&0&0&0&-4&0\cr
\noalign{\vskip5pt\hrule\vskip4pt}
\multispan{12}\hfil\hbox{Table 11}\hfil\cr
}}$$

\th EXERCISE 12
\enonce
Let\/ $n$ be an odd squarefree integer.  If\/ $n\equiv5,7\pmod8$, then\/
$c_1(n)=0$.  If\/ $n\equiv3\pmod4$, then\/ $c_2(n)=0$.
\endth

{\it For the remainder of this section, let $\alpha$ be a squarefree 
  integer $>0$, and write $\alpha=jn$, with $j=1,2$ and $n$ odd.}

\th CONJECTURE 13
\enonce
The number\/ $\alpha=jn$ is congruent if and only if\/ $c_j(n)=0$. 
\endth

As we shall see, this conjecture is implied by Conjecture~24 (Birch and
Swinnerton-Dyer), combined with Theorem~25 (Tunnell).

The reader should marvel at how unexpected the (conjectural) characterisation
is, how far-removed from rational right triangles and their areas~!  

{\eightpoint The physicist Richard Feynman claims in his {\it Surely you are
    joking\/} that he could guess whether a mathematical statement explained
  to him in elementary terms was true or false.  It would have been
  interesting to have given him Definition~1 and Notation~10, and to have
  asked him if Conjecture~13 is true.}

We do know one of the implications in Conjecture~13~:

\th THEOREM 14 (J. Coates \& A. Wiles \citer\coateswiles())
\enonce
If $c_j(n)\neq0$, then $\alpha=jn$ is not congruent.
\endth

It follows for example that the numbers $1, 2, 3, 10, 17, 19, 26$ and $42$
(Table~11) are not congruent (cf.~Theorem~3).

If the squarefree odd integer $n$ is $\equiv3\pmod4$
(resp.~$\equiv5,7\pmod8$), then $2n$ (resp.~$n$) should be congruent
(Exercise~12, Conjecture~13), and the first few such $n$ are indeed so
(cf.~(3)).  In a paper which became influencial when it was properly
understood, an obscure schoolteacher by the name of K.~Heegner proved that
this is true if $n$ is prime \citer\heegner().  

However, in general, the result is only conditional.  It is conditional on the
finiteness of a certain set $S_\alpha$, which will be discussed in detail in
later parts of this report (cf.~Conjecture~16).  Suffice it to say here that
the finiteness of the set $S_\alpha$ is equivalent to the finiteness of the
group ${\cyr Sh}(E_\alpha)$ which is more familiar to arithmeticians.  We have
chosen to formulate things in terms of $S_\alpha$, which can be defined in
elementary terms.

\th THEOREM 15 (K. Rubin \citer\rubin())
\enonce
If\/ $c_j(n)=0$ and if the set\/ $S_{\alpha}$ is finite, then the number\/
$\alpha=jn$ is congruent.  
\endth

%{\eightpoint A detailed discussion of the sets $S_{\alpha}$ is given below.}

Note that if $c_j(n)=0$ and if $S_{\alpha}$ is finite, then Theorem~15 shows
that $\alpha$ is congruent {\it without\/} exhibiting a rational right
triangle of area $\alpha$.  However, in some cases (``rank 1''), the set
$S_\alpha$ is known to be finite and there is a method (``Heegner points'')
for constructing such a triangle.  Zagier's example showing that 157 is
congruent, displayed before Theorem~3, is of this type.  

\bigbreak
\centerline{\bf 3. Local number fields}
\medskip

The group ${\bf Q}^\times$ (modulo its torsion subgroup $\{1,-1\}$) admits the
set of prime numbers as a ${\bf Z}$-basis.  For every prime number $p$, there
is thus a unique homomorphism $v_p:{\bf Q}^\times\rightarrow{\bf Z}$ such that
$v_p(p)=1$ and $v_p(l)=0$ for every prime number $l\neq p$~; extending it to
${\bf Q}$ by $v_p(0)=+\infty$, we get a {\it discrete valuation}, because it
satisfies 
$$
v_p(x+y)\ge\inf(v_p(x),v_p(y))\qquad\hbox{for all\ } x,y\in{\bf Q}.
$$
Define $|\ |_p:{\bf Q}\rightarrow{\bf R}$ by $|x|_p=p^{-v_p(x)}$
(convention~: $p^{-\infty}=0$).  Then $|x-y|_p$ is a {\it distance\/} on
${\bf Q}$ with respect to which it can be completed to obtain a field ${\bf
  Q}_p$ much in the same way as we obtain the field ${\bf R}$ from ${\bf Q}$
by completing it with respect to the usual distance
$|x-y|_\infty=\sup(x-y,y-x)$.  For this reason, the field of real numbers is
sometimes denoted ${\bf Q}_\infty$.

{\eightpoint It can be shown that the $v_p$ ($p$ prime) are the only discrete
  valuations, and $|\ |_\infty$ the only archimedean absolute value, on the
  field ${\bf Q}$ (A. Ostrowski, 1918).  Thus the absolute values $|\ |_p$
  ($p$ prime or $p=\infty$) determine all the {\it places\/} of $\bf Q$.}

The fields ${\bf Q}_p$ (including $p=\infty$) play a fundamental role in
arithmetic.  It is always a good idea to first study ``global questions'' ---
questions about rational numbers --- everywhere ``locally'' in the fields
${\bf Q}_p$, before trying to answer the original question.  We discuss a
basic example in the next section.

For $p$ prime, the field ${\bf Q}_p$ comes equipped with a continuous discrete
valuation extending $v_p$; elements of positive valuation form a subring ${\bf
  Z}_p$ (``the ring of integers'') in which $p{\bf Z}_p$ is the unique maximal
ideal.  The quotient ${\bf Z}_p/p{\bf Z}_p$ is the finite field ${\bf F}_p$
(``the residue field'') of $p$ elements.

{\eightpoint All books on Number Theory (Artin, Hasse, Weil, Serre,
  Kato-Kurokawa-Saito, $\ldots$) provide an introduction to the fields ${\bf
    Q}_p$ and their extensions.}

\bigbreak 
\centerline{\bf 4. The local-to-global principle for conics} 
\medskip

To avoid speaking of curves, we use the equivalent language of a {\it function
  field\/} $F$ over a field $k$~: a finitely generated extension of $k$ in
which $k$ is algebraically closed~; we'll be mostly concerned with the case
when $F$ has transcendence degree~$1$ over~$k$.  Concretely, if $f\in k[x,y]$
is an absolutely irreducible polynomial --- one which remains irreducible over
every finite extension of $k$ ---, then the field of fractions $F$ of the
(integral) ring $k[x,y]/fk[x,y]$ is a function field over $k$~; we write
$F=k(x,y)$, with the relation $f=0$.  For every extension $L$ of $k$, we then
get a function field over $L$ by ``extending the scalars'' of $F$ from $k$ to
$L$~: the field of fractions of $L[x,y]/fL[x,y]$.

Let us fix an algebraic closure $\bar{\bf Q}$ of $\bf Q$.  Clearly, the
function field ${\bf Q}(T)$ becomes isomorphic to $\bar{\bf Q}(T)$ over
$\bar{\bf Q}$.  Are there any others which do~?  And, is there a way to
classify them all~?

Fix an algebraic closure $\bar{\bf Q}_p$ of ${\bf Q}_p$.  The corresponding
local question is~: find all function fields over ${\bf Q}_p$ which become
isomorphic over $\bar{\bf Q}_p$ to $\bar{\bf Q}_p(T)$.  Such functions fields
will be called {\it solutions\/} to our problem.

The trivial solution to the problem is the rational function field ${\bf
  Q}_p(T)$.  {\it It can be shown that there is precisely one other
  solution}~; let us call it $F_p$.  Thus the function field $F_p$ is not the
rational function field but becomes (isomorphic to) the rational function
field over $\bar{\bf Q}_p$.  For example, when $p=\infty$, the field
$F_\infty$ is ${\bf Q}_\infty(x,y)$ with the relation $x^2+y^2+1=0$.
Moreover, for every place $p$, it is an easy matter to decide if a given
``local solution'' is isomorphic ${\bf Q}_p(T)$ or to $F_p$.
 
Now, if $F$ is ``global solution'' to our problem, then it is a ``local
solution'' everywhere.  In other words, if $F$ is a function field over ${\bf
  Q}$ which becomes the rational function field over $\bar{\bf Q}$, then $F$
becomes isomorphic to one of ${\bf Q}_p(T)$, $F_p$ over every completion ${\bf
  Q}_p$ of ${\bf Q}$, including $p=\infty$.

{\it What can be shown is that every\/ $F$ becomes isomorphic to\/ ${\bf
    Q}_p(T)$ for almost every\/ $p$, the places where it doesn't --- there are
  thus only finitely many of them --- are even in number, and, given any
  finite set\/ $\Sigma$ of places, even in number, there is a unique global
  solution which becomes isomorphic to\/ $F_p$ for all\/ $p\in\Sigma$ and to\/
  ${\bf Q}_p(T)$ for all\/ $p\notin\Sigma$.}

\goodbreak

{\eightpoint There are many equivalent ways --- curves of genus~0, quadratic
  forms in three variables, quaternion algebras --- of expressing this
  principle.}

It follows that if two global solutions $F$, $F'$ are ``everywhere locally
isomorphic'' (become isomorphic to each other at every place $p$, including
$p=\infty$), then they are ${\bf Q}$-isomorphic.  This happy circumstance is
expressed by saying that such function fields obey {\it the local-to-global
  principle}.  (In fact, in the case at hand, it is sufficient to demand that
$F$, $F'$ be isomorphic at all places but one~; they are then automatically
isomorphic at the remaining place.)

{\eightpoint The best accounts of this circle of ideas, in the equivalent
  language of quadratic forms, are to be found in Serre's {\it Course in
    arithmetic\/} and in {\it Number theory 1, Fermat's dream\/} by Kato,
  Kurokawa and Saito.  A theorem of Adrien-Marie Legendre can be considered to
  be a precursor of local-to-global considerations, see Weil's {\it Number
    theory, an approach through history.}}
  
% The $2$-variable function field $F$ defined by $y^2+3z^2=x^3-2$ is not
% isomorphic to ${\bf Q}(T_1,T_2)$ but becomes isomorphic to ${\bf
%   Q}_p(T_1,T_2)$ at every place $p$ because $F(T_3,T_4,T_5)$ is the
% $5$-variable rational function field over ${\bf Q}$ (A. Beauville, J.-L.
% Colliot-Th{\'e}l{\`e}ne, J.-J. Sansuc, H. P. F. Swinnerton-Dyer).  N.
% Shepherd-Barron has recently shown that adjoining $T_5$ is superfluous.

% Not true. Cf. the exchange of emails with Colliot, 27/2/2007.

{\eightpoint I don't know of any classification of function fields over ${\bf
    Q}$ which become the $2$-variable rational function field over every
  completion.}

\bigbreak
\centerline{\bf 5. The failure of the local-to-global principle}
\medskip

In the last section we saw that the local-to-global principle holds for
functions fields over ${\bf Q}$ which become isomorphic over $\bar{\bf Q}$ to
the rational function field.  Such function fields are of the form ${\bf
  Q}(x,y)$, $ax^2+by^2=1$, for some $a,b\in{\bf Q}^\times$, and it is easy to
decide when this field is isomorphic to the one defined by $a'x^2+b'y^2=1$
($a',b'\in{\bf Q}^\times$), because it suffices to check that they are
isomorphic everywhere locally.

In the early 1940s, Carl-Erik Lind and Hans Reichardt found the first examples
of function fields which violate the local-to-global principle.  Equivalently,
Reichardt showed that $2y^2=1-17x^4$ has solutions in every completion of
${\bf Q}$ but no rational solutions --- not even ``at infinity''.

{\eightpoint Lind's thesis was reviewed by Andr{\'e} Weil in the {\it
    Mathematical Reviews}, and it is amazing to note that he does not mention
  this discovery.  Nor does the reviewer of Reichardt's paper, in spite of the
  explicit title~: {\it Einige im Kleinen {\"u}berall l{\"o}sbare, im Grossen
    unl{\"o}sbare diophantische Gleichungen}.  It must be said that the first
  instance of the failure of a local-to-global principle, due to Hasse, was
  discovered by him after he had proved its validity for quadratic forms.}

The example most commonly cited these days, originating with Ernst Selmer,
is that of the function field ${\bf Q}(x,y)$, $3x^3+4y^3+5$.  Cf.~Example~27. 

Let $\alpha$ be a squarefree integer $>0$ and consider the function field
${\bf Q}(C_\alpha)$ defined by the equation $C_\alpha:\alpha y^2=x^3-x$.  It
may happen that there are many function fields $F$ over ${\bf Q}$ which become
isomorphic to ${\bf Q}(C_\alpha)$ at every place $p$ of ${\bf Q}$.  In other
words, ${\bf Q}(C_\alpha)$ may have ``twisted forms'' $F$ which become
isomorphic to it when we extend scalars of $F$ and ${\bf Q}(C_\alpha)$ from
${\bf Q}$ to ${\bf Q}_p$.  Let us denote the set of isomorphism classes of
such $F$ by $S_\alpha$.  This is the set which appears in Theorem~15.

Thus the problem of congruent numbers would be solved if we could settle the
following conjecture, whose generalisation Conjecture~26 is a major open
question in contemporary arithmetic.  

\th CONJECTURE 16 (I. Shafarevich \& J. Tate)
\enonce
For every\/ $\alpha$, the set\/ $S_\alpha$ of function fields which become
isomorphic to\/ ${\bf Q}(C_\alpha)$ at every place is finite.
\endth

The more standard version of this conjecture asserts the finiteness of the
{\it group\/} ${\cyr Sh}(E_\alpha)$, whose definition is more advanced.  The
reader who knows it should be able to prove that $S_\alpha$ is finite if and
only if ${\cyr Sh}(E_\alpha)$ is finite \citer\mazur().  The same remark
applies to Conjecture~26.  

We have seen that the congruent number problem amounts to the arithmetic study
of the equation $\alpha y^2=x^3-x$, which can be rewritten as
$y^2=x^3-\alpha^2 x$.  The rest of this report is devoted to a rapid survey of
the arithmetic of equations of the type $y^2=f(x)$, where $f\in{\bf Q}[x]$ is
a monic cubic polymonial with distinct roots (in $\bar{\bf Q}$).

\bigbreak
\centerline{\bf 6. Elliptic curves : results and conjectures}
\medskip

In the next two sections, we enumerate some arithmetic properties of elliptic
curves.  For the sake of simplicity, we work over the field ${\bf Q}$~; the
only exceptions being a result over finite fields, one over ${\bf Q}_p$, and
an example over ${\bf Q}(i)$.

An {\it elliptic curve\/} $E$ over a field $k$ is a curve defined in the
projective plane by an equation of the type
$$
f(x,y)=y^2+a_1xy+a_3y-x^3-a_2x^2-a_4x-a_6=0\quad (a_i\in k)\leqno{(5)}
$$
without singularities, a condition which says that the discriminant
$\Delta$ --- a certain polynomial in the $a_i$ --- is $\neq0$, or equivalently
that the corresponding function field is of ``genus~$1$'', unlike the function
fields which become isomorphic to $\bar{\bf Q}(T)$, which are of genus~0.  

More precisely, the discriminant of $f$ --- the result of elliminating
$x$, $y$ from $f$, $f'_x$, $f'_y$ --- is 
$$
\Delta=-b_2^2b_8-2^3b_4^3-3^3b_6^2+3^2b_2b_4b_6
$$
where
$$
b_2=a_1^2+2^2a_2,\quad b_4=a_1a_3+2a_4,\quad b_6=a_3^2+2^2a_6
$$
and
$$
b_8=b_2a_6-a_1a_3a_4+a_2a_3^2-a_4^2.
$$

The curve $E$ has a ``point at infinity'' $O$~; for any extension $L$ of $k$,
there is a natural group law on the set $E(L)$ consisting of $O$ and the
solutions of (5) in $L$, uniquely determined by the requirement that $O$ be
the origin and that the sum of the three points (counted with multiplicity) in
which $E$ intersects a given line be $O$~; the groups $E(L)$ are commutative.
Elements of $E(L)$ can be identified with triples $(x,y,z)\neq(0,0,0)$
($x,y,z\in L$) satisfying the homogenised version of $(5)$; two such triples
being considered the same if each is a multiple of the other by an element of
$L^\times$.  Two elliptic curves are isomorphic if the corresponding function
fields over $k$ are $k$-isomorphic.

{\eightpoint Let $C$ be a smooth proper absolutely connected genus-$1$
  $k$-curve and let $J$ be its jacobian --- a $1$-dimensional abelian
  $k$-variety.  If $C$ has a $k$-rational point $O$, there is a unique
  $k$-morphism $C\rightarrow J$ sending a point $P$ to the class of the
  divisor $P-O$~; it is an isomorphism.}

{\eightpoint For surveys of arithmetic on elliptic curves, see Cassels
  \citer\cassels() and Tate \citer\tate().}

\th THEOREM 17 (L. Mordell, 1922)
\enonce
For every elliptic curve\/ $E$ over\/ ${\bf Q}$, the group\/ $E({\bf Q})$ is 
finitely generated.
\endth

This result was conjectured by Henri Poincar{\'e} around 1900.  Mordell's
proof is a generalisation of Fermat's method of infinite descent --- employed
in the proof of Theorem~3 ---~; its modern renditions consist of two parts.

The first part shows that the group $E({\bf Q})/2E({\bf Q})$ is finite.  The
second part studies a canonical real-valued ``height'' function $h$ on $E({\bf
  Q})$, coming from the various absolute values of ${\bf Q}$.  The method of
infinite descent is distilled in the statement that a commutative group
$\Gamma$, endowed with such a function $h$, for which $\Gamma/2\Gamma$ is
finite, is necessarily finitely generated.
  
{\eightpoint Accounts of the proof can be found in the books by Weil and
  Kato-Kurokawa-Saito cited above, as well as in Silverman-Tate, {\it Rational
    points on elliptic curves}.}

By contrast, the corresponding local result says that for an elliptic curve
$E$ over $\bf R$, the group $E({\bf R})$ has a subgroup of index at most $2$
isomorphic to ${\bf R}/{\bf Z}$, and, for an elliptic curve $E$ over ${\bf
  Q}_p$ ($p$ prime), $E({\bf Q}_p)$ has a subgroup of finite index isomorphic
to ${\bf Z}_p$.  Thus, for an elliptic curve $E$ over ${\bf Q}$, although the
three groups ${\bf Q}$, ${\bf Q}^\times$, $E({\bf Q})$ have very different
structures, they are ``almost the same'' everywhere locally.  However, the
indices in question are important local invariants of $E$.

For a given $E$ over ${\bf Q}$, the torsion subgroup of $E({\bf Q})$ is easy
to determine ({\'E}. Lutz)~; for example, the torsion subgroup of
$C_\alpha({\bf Q})$ consists of $O$ and the three points $(-1,0)$,
$(0,0)$, $(1,0)$ of order~$2$.

\th THEOREM 18 (B. Mazur \citer\eisenstein())
\enonce
Let\/ $E$ be an elliptic curve over\/ ${\bf Q}$.  The torsion subgroup of\/
$E({\bf Q})$ is isomorphic to one of the fifteen groups
$$
{\bf Z}/m{\bf Z}\ (m=1,2,\ldots,10,12),\quad
{\bf Z}/2{\bf Z}\times{\bf Z}/2\nu{\bf Z}\ (\nu=1,2,3,4).
$$
% and each of these groups occurs infinitely often.
\endth

No {\it uncondtional\/} method is known, however, for determining the rank of
$E({\bf Q})$ for a given $E$.  The set of possible ranks for variable $E$ is
not known either, but N.~Elkies has recently produced examples where
$\mathop{\rm rk}E({\bf Q})$ is at least $28$.  We shall mostly concentrate on
the question of deciding if the rank is $0$ or $>0$.

Let $p$ be a prime number and let $E$ be an elliptic curve over ${\bf Q}_p$,
given by an equation $(5)$.  We may assume by a change of variables that
$a_i\in{\bf Z}_p$~; the discriminant $\Delta$ is then in ${\bf Z}_p$.  If the
$a_i\in{\bf Z}_p$ can be so chosen that $\Delta\in{\bf Z}_p^\times$, we say
that $E$ has {\it good reduction\/} at $p$~; if so, the equation $(5)$, read
modulo~$p$, defines an elliptic curve $E_p$ --- uniquely determined by $E$ and
$p$ --- over the finite field ${\bf F}_p$, and there is a homomorphism $E({\bf
  Q}_p)\rightarrow E_p({\bf F}_p)$ which sends a point to the reduction
modulo~$p$ of any of its representatives $(x,y,z)$ with coordinates in ${\bf
  Z}_p$ and at least one coordinate in ${\bf Z}_p^\times$.

There is a criterion for good reduction (``N{\'e}ron-Ogg-Shafarevich'').  Let
$\bar{\bf Q}_p$ be an algebraic closure of ${\bf Q}_p$.  There is a unique
extension of $v_p$ to a valuation $v_p:\bar{\bf Q}_p^\times\rightarrow{\bf Q}$
of which the residue field $\bar{\bf F}_p$ is an algebraic closure of ${\bf
  F}_p$.  The {\it inertia group\/} is the kernel of the natural surjection
$\mathop{\rm Gal}(\bar{\bf Q}_p|{\bf Q}_p)\rightarrow\mathop{\rm Gal}(\bar{\bf
  F}_p|{\bf F}_p)$~; it acts on the $m$-torsion\/ ${}_mE(\bar{\bf Q}_p)$ for
every $m$.
 
\th THEOREM 19 (J.-P. Serre \& J. Tate \citer\serretate())
\enonce
An elliptic curve\/ $E$ over\/ ${\bf Q}_p$ has good reduction if and only if
the the action of the inertia group on\/ ${}_mE(\bar{\bf Q}_p)$ is trivial
for every $m$ prime to $p$. 
\endth

Every elliptic curve $E$ over ${\bf Q}$ has good reduction at almost all
primes.  One might ask to what extent $E$ is determined by the the number
$|E_p({\bf F}_p)|$ of points modulo $p$ for varying $p$.  We say that two
elliptic curves are {\it isogenous\/} if their function fields can be embedded
into each other.

\th THEOREM 20 (G. Faltings \citer\faltings())
\enonce
If\/ $E'$ is an elliptic curve over\/ ${\bf Q}$ such that\/  $|E'_p({\bf
  F}_p)|=|E_p({\bf F}_p)|$ for almost all primes\/ $p$, then\/ $E'$ is
isogneous to\/~$E$.  
\endth

There are only finitely many such $E'$ because isogenous curves are known to
have good reduction at the same primes (cf.~Theorem~19)~:

\th THEOREM 21 (I. Shafarevich, 1962)
\enonce
Given a finite set\/ $T$  of primes, there are only finitely many elliptic 
curves over\/ ${\bf Q}$ having good reduction at every prime\/ $p\notin T$.
\endth

There is a sense in which the more fundamental quantity is not $|E_p({\bf
  F}_p)|$ but $a_p(E)$, defined by $|E_p({\bf F}_p)|=1-a_p(E)+p$, and there is
sense in which the following theorem is the analogue, for function fields of
elliptic curves over finite fields, of the famous Riemann Hypothesis~: ``the
zeros in the critical strip $0<\mathop{\rm Re}(s)<1$ of the zeta function
$\zeta$  of ${\bf Q}$ have real part $1\over 2$''.

\th THEOREM 22 (H. Hasse, 1933)
\enonce
Let\/ $A$ be an elliptic curve over a finite field\/ $k$ of\/ $q$ elements.
Define the integer\/ $a$ by $|A(k)|=1-a+q$.  Then $|a|\le 2\sqrt q$.
\endth

Returing to our $E$ over ${\bf Q}$, Birch and Swinnerton-Dyer argued that if\/
$E({\bf Q})$ is infinite, the groups $E_p({\bf F}_p)$ (for $p$ a prime of good
reduction for $E$) should have more elements ``on the average'' than if
$E({\bf Q})$ is finite.  In view of Hasse's theorem, the product $\prod_p
{p\over |E_p({\bf F}_p)|}$ should diverge to $0$ if the rank is $>0$, and
converge to a limit $\neq0$ if the rank is $0$.  This is made precise in terms
of the $L$-function of $E$.
  
For a prime $p$ of good reduction for $E$, we have the number $a_p(E)$~; for
``cohomological'' reasons, consider the infinite product (for $s\in{\bf C}$) 
$$
L(E,s)=\prod_{p}{1\over 1-a_p(E).p^{-s}+p.p^{-2s}}\raise2pt\hbox{.}
$$
Theorem 22 implies that this converges for $\mathop{\rm Re}(s)>{3\over2}$, but
more is true~:

\th THEOREM 23 (A. Wiles, R. Taylor, F. Diamond, B. Conrad, C. Breuil, 
1995--2000) 
\enonce
The function\/ $L(E,s)$ admits an analytic continuation to the whole of\/
${\bf C}$.
\endth

{\eightpoint For the congruent number elliptic curves $C_\alpha$, this is due
  to Andr{\'e} Weil.  There is a way of introducing factors in $L(E,s)$
  corresponding to the primes which divide $\Delta$, and indeed to the place
  $\infty$.  This ``completed'' $L$-function $\Lambda(E,s)$ has a ``functional
  equation'' for $s\mapsto 2-s$, just as $\zeta$ has a functional equation for
  $s\mapsto s-1$.}

Note that the product $\prod_p {p\over |E_p({\bf F}_p)|}$ is formally equal to
$L(E,1)$, and the above heuristic considerations and extensive calculations on
one of the first electronic computers at Cambridge led to the following
conjecture.

\th CONJECTURE 24 (B. Birch \& P. Swinnerton-Dyer, 1965)
\enonce
The group\/ $E({\bf Q})$ is infinite if and only if\/ $L(E,1)=0$.  More
precisely, its rank equals the order of vanishing of\/ $L(E,s)$ at\/ $s=1$. 
\endth

{\eightpoint The order of vanishing of the completed $L$-function
  $\Lambda(E,s)$ is the same as that of $L(E,s)$ at $s=1$.  There is a refined
  version of Conjecture~24 which gives the leading coefficient of
  $\Lambda(E,s)$ at $s=1$ in terms of the local and global arithmetic
  invariants of the curve $E$~; its formulation is subject to the truth of
  Conjecture~26.}

Conjecture 13 follows from this, thanks to the following criterion~: 

\th THEOREM 25 (J. Tunnell \citer\tunnell())
\enonce
For a squarefree integer\/ $\alpha=jn$ ($j=1,2$ and\/ $n$ odd), one has\/
$L(C_\alpha,1)=0$ if and only if\/ $c_j(n)=0$.
\endth

The elliptic curve $E$ has the function field ${\bf Q}_l(E)$ at the various
places $l$ of ${\bf Q}$.  Just as we did in the case of the congruent number
elliptic curves $C_\alpha$, we now consider the set $S_E$ of (isomorphism
classes of\/) all function fields over ${\bf Q}$ which becomes isomorphic to
${\bf Q}_l(E)$ at every place $l$~; of course, ${\bf Q}(E)$ belongs to $S_E$.

\th CONJECTURE 26 (I. Shafarevich \& J. Tate)
\enonce
For every elliptic curve\/ $E$ over\/ ${\bf Q}$, the set\/ $S_E$ is finite.
\endth

The original conjecture asserts the finiteness, for every $E$ over ${\bf Q}$,
of the {\it group\/} ${\cyr Sh}(E)$ of ``torsors'' under $E$ which are
``everywhere locally trivial''.  This is equivalent to the finiteness of
$S_E$.

% {\eightpoint The correct language to study the set $S_\alpha$ is that of
%   curves, their automorphisms, and galoisian cohomology.}

{\eightpoint Yuri Manin has introduced an ``obstruction'' to explain the
  failure of the local-to-global principle for the function field ${\bf Q}(E)$
  of an elliptic curve $E$ over ${\bf Q}$.  He shows that the finiteness of
  $S_E$ is equivalent to his obstruction being the only one.}

The equation $x^3+y^3+60=0$ can be put in the form $(5)$ by a change of
variables~; it therefore defines an elliptic curve.  

\th EXAMPLE 27 (Mazur \citer\mazur())
\enonce
For\/ $E$ defined by\/ $x^3+y^3+60=0$, the set\/ $S_E$ consists
of\/ ${\bf Q}(E)$ and the function fields
$$
3x^3+4y^3+5,\ 12x^3+y^3+5,\ 15x^3+4y^3+1,\ 3x^3+20y^3+1.
$$
\endth

The best available result in the direction of Conjectures~24 and~26 to date,
the fruit of a succession of papers by numerous mathematicians, is a theorem
of Victor Kolyvagin, of which the theorem of Coates \& Wiles (Theorem~14) is a
particular case, and which subsumes some of the results of Benedict Gross and
Zagier \citer\grosszagier().

\th THEOREM 28 (V. Kolyvagin \citer\kolyvagin())
\enonce
If\/ $L(E,1)\neq0$, then\/ $E({\bf Q})$ is finite.  If\/ $L(E,s)$ has a simple
zero at\/ $s=1$, then\/ $E({\bf Q})$ has rank\/ $1$.  In both these cases, the
set\/ $S_E$ is finite.
\endth

If the zero at $s=1$ has multiplicity $>1$, Conjecture~26 is needed
(Cf. Theorem~15)~: 

\th THEOREM 29 (C. Skinner \& {\'E}. Urban \citer\skinnerurban())
\enonce
Suppose that\/ $L(E,1)=0$ and that the set\/ $S_E$ is finite.  Then the
group\/ $E({\bf Q})$ is infinite.
\endth

{\eightpoint There is a parallel theory of elliptic curves $E$ over function
  fields $F$ over finite fields.  The analogue of Mordell's theorem
  (Theorem~17) is true~: the group $E(F)$ is finitely generated.  K. Kato \&
  F. Trihan \citer\katotrihan() have proved the analogue of (the refined
  version of) the Birch \& Swinnerton-Dyer conjecture (Conjecture~24), subject
  to the truth of the analogue of the Shafarevich-Tate conjecture
  (Conjecture~26).}

The study of ``special values'' of $L$-functions, of which the refined
conjecture of Birch and Swinnerton-Dyer is the prototype, is one of the major
themes of contemporary arithmetic.

\bigbreak
\centerline{\bf 7. Complex multiplications}
\medskip

Let $E$ be an elliptic curve over $\bar{\bf Q}$.  Because $E$ has a group law,
there are many embeddings of the function field $\bar{\bf Q}(E)=\bar{\bf
  Q}(x,y)$ into itself~: for every integer $n\neq0$, there is an embedding
$[n]_E$ which sends $x,y$ to $x_n,y_n$, the coordinates of the multiple $nP$
of the point $P=(x,y)$~; it is of degree $n^2$.  For example, when $E$ is the
congruent number elliptic curve $C_\alpha$~$(4)$ and $n=-1$, it is the
automorphism $x\mapsto x$, $y\mapsto-y$ of the function field~; when $n=-2$,
it is the degree-$4$ embedding given in Exercise~7.

In a sense, for most elliptic curves, these are the only embeddings of the
function field into itself.  But there are some elliptic curves for which
there are more embeddings, for example the automorphism $x\mapsto-x$,
$y\mapsto iy$ ($i$ being a chosen square root of $-1$) of $\bar{\bf
  Q}(C_\alpha)$ whose square is $[-1]_{C_\alpha}$.  In such a case we say that
the elliptic curve $E$ has ``complex multiplications''~; it then determines an
imaginary quadratic field $K$, which in the case of the $C_\alpha$ is ${\bf
  Q}(i)$~; we say that $E$ has complex multiplications by $K$.

The arithmetic properties of elliptic curves differ vastly according as they
have complex multiplications or not.  For example, for elliptic curves having
complex multiplications, the theorem ``$L(E,1)\neq0\Rightarrow E({\bf
  Q})\hbox{ is finite}$'' was proved for by Coates-Wiles (cf.~Theorem~14) a
good eleven years before Kolyvagin's general result (cf.~Theorem~28), the
analytic continuation of $L(E,s)$ was proved by Weil and Max Deuring in
1953--1957, much before the general result of Wiles and his school in
1995--2000 (Theorem~23), and the implication ``$L(E,1)=0\hbox{ and } S_E\hbox{
  finite }\Rightarrow E({\bf Q})\hbox{ infinite}$'' was proved by Rubin
(cf.~Theorem~15) some fifteen years before the general result of Skinner-Urban
(Theorem~29).

We illustrate the differences by three examples.  For the first, recall that
an elliptic curve $A$ over ${\bf F}_p$ is called {\it supersingular\/} if the
$p$-torsion ${}_pA(\bar{\bf F}_p)$ is reduced to $\{O\}$, or, equivalently for
$p\neq2,3$, if $|A({\bf F}_p)|=1+p$ (equivalently, $a=0$, in the notation of
Theorem~22).  Returning to our $E$ over ${\bf Q}$, we ask~: How often is $E_p$
supersingular~?  Deuring showed if $E$ has complex multiplications, then this
happens for half the primes $p$ (cf.~Example~33); if not, Jean-Pierre Serre
proved that the set of primes in question has density~$0$.  That it is
infinite is a relatively recent result.

\th THEOREM 30 (N. Elkies \citer\elkies())
\enonce
For every elliptic curve\/ $E$ over\/ ${\bf Q}$, there are infinitely many
primes\/ $p$ at which $E_p$ is supersingular.
\endth

For the second example, recall that for every prime $p$, if we adjoin the
$p$-torsion of the multiplicative group $\bar{\bf Q}^\times$, which consists
of $p^{\rm th}$ roots of~1, to ${\bf Q}$, we get a galoisian extension ${\bf
  Q}({}_p\mu)$ whose group of automorphisms is $\mathop{\rm Gal}({\bf
  Q}({}_p\mu)|{\bf Q})=GL_1({\bf F}_p)$.  For an elliptic curve $E$ over ${\bf
  Q}$, the $p$-torsion of $E(\bar{\bf Q})$ is a $2$-dimensional vector ${\bf
  F}_p$-space~; if we adjoin it to ${\bf Q}$, we get a galoisian extension
${\bf Q}({}_pE)$.  What is $\mathop{\rm Gal}({\bf Q}({}_pE)|{\bf Q})$~?

\th THEOREM 31 (J.-P. Serre \citer\serre())
\enonce
Suppose that\/ $E$ does not have complex multiplications.  Then the group of 
automorphisms of\/ ${\bf Q}({}_pE)$ is\/ $GL_2({\bf F}_p)$ for almost all ---
all but finitely many --- primes\/ $p$. 
\endth

The corresponding local result for $E$ over ${\bf Q}_l$ says, at least in the
case of good reduction, that $\mathop{\rm Gal}({\bf Q}_l({}_pE)|{\bf Q}_l)$ is
cyclic for $l\neq p$ (cf.~Theorem~19).

If $E$ (over ${\bf Q}$) has complex multiplications, the group of
automorphisms is much smaller~: if $K$ --- an imaginary quadratic field --- is
the field of complex multiplications, then $K({}_pE)$ is an abelian extension
of $K$.  However, such $E$ serve a different, if related, purpose.

Recall that the theorem of Kronecker-Weber asserts that {\it if we adjoin the
  entire torsion subgroup of\/ $\bar{\bf Q}^\times$ --- all roots of\/ $1$ ---
  to\/ ${\bf Q}$, we get the maximal abelian extension.}  Generating the
maximal abelian extension of other number fields is a major open problem
(Kronecker's {\it Jugendtraum}, Hilbert's Problem~12)~; the theory of complex
multiplications provides the answer in the case of imaginary quadratic fields,
as in the next example.

\th EXAMPLE 32
\enonce
Let\/ $E$ be the elliptic curve\/ $y^2=x^3+x$, which has complex
multiplications by\/ ${\bf Q}(i)$.  If we adjoin the entire torsion subgroup
of\/ $E(\bar{\bf Q})$ to ${\bf Q}(i)$, we get the maximal abelian extension
of\/ ${\bf Q}(i)$.
\endth

Our third example concerns a ``formula'' for $a_p(E)$ for a fixed $E$ and
varying $p$.  There is indeed such a formula if $E$ has complex
multiplications, as illustrated by a theorem of Carl Gauss about the curve
$x^3+y^3+1=0$ (which can be put in the canonical form $(5)$, and has complex
multiplications by ${\bf Q}(j)$, $j^2+j+1=0$).  It uses the fact that for a
prime $p\equiv1\pmod3$, there is a pair of integers $(c_p,d_p)$, unique up to
signs, such that $4p=c_p^2+27d_p^2$~; to fix the sign of $c_p$, assume that
$c_p\equiv-1\pmod3$.
 
\th EXAMPLE 33 (C. Gauss, 1801) 
\enonce 
Let\/ $E$ be the elliptic curve\/ $x^3+y^3+1=0$ and $p$ a prime.  If\/
$p\equiv1\pmod3$, then\/ $a_p(E)=c_p$. If\/ $p\equiv-1\pmod3$, then\/
$a_p(E)=0$.  
\endth
 
{\eightpoint See Silverman-Tate for a proof.  Note that this implies
  Theorem~22 for $E$.}

By contrast, if $E$ does not have complex multiplications, the behaviour of
the $a_p(E)$ is entirely different.  Mikio Sato and Tate independently arrived
at a conjectural distribution law for $\gamma_p(E)=a_p(E)/2\sqrt p$, which
lies between $-1$ and $+1$ for every $p$ (cf.~Theorem~22).  How often does it
lie in $[\beta,\delta]\subset[-1,+1]$~?

\th CONJECTURE 34 (M. Sato \& J. Tate, 1960)
\enonce
Suppose that\/ $E$ does not have complex multiplications, and let\/
$[\beta,\delta]\subset[-1,+1]$ be an interval.  Then the proportion of
primes\/ $p$ for which\/ $\gamma_p(E)\in[\beta,\delta]$ is given by
$$
{2\over\pi}\!\int_{\!\beta}^\delta\!\!\!\sqrt{1-x^2}\;dx.
$$
\endth
This conjecture has been proved, subject to a mild technical hypothesis on
$E$, by Laurent Clozel, Michael Harris, Nicholas Shepherd-Barron and Richard
Taylor in a series of three papers in early~2006.  The technical hypothesis
demands that $E$ have ``multiplicative reduction'' at some prime~$p$, which
means roughly that the best possible reduction at $p$ is not an elliptic curve
$E_p$ as in the case of good reduction, but the multiplicative group (and not
the additive group --- the third possibility).  An algorithm due to Tate
allows one to determine the type of reduction at any given~$p$ in terms of the
coefficients $a_i$ $(5)$ defining $E$.  Concretely, although we cannot choose
$a_i\in{\bf Z}_p$ with minimal $v_p(\Delta)$ so as to have $v_p(\Delta)=0$,
they can be so chosen as to have $v_p(c_4)=0$, where $c_4=b_2^2-2^3.3.b_4$, and
the $b_i$ are displayed after equation~$(5)$.  It is only a matter of time
before this hypothesis is removed.

\th THEOREM 35 (L. Clozel, M. Harris, N. Shepherd-Barron \& R. Taylor
\citer\taylor())
\enonce
Conjecture 34 is true if\/ $E$ has multiplicative reduction at some
prime\/~$p$. 
\endth

\bigbreak
\centerline{\bf 8. Modular forms}
\medskip

We have not mentioned them, although they have appeared here without being
named.  If we evoke them here, it is only to say that most of the spectacular
recent results which we have enumerated would not have been possible without
their help.  Take the analytic continuation of $L(E,s)$ (Theorem~23)~: the
crucial result (Wiles and others) is to show that the sequence $(a_p(E))_p$
defines a modular form.

Results of Gross-Zagier and of Kolyvagin (Theorem~28), which predate Wiles,
were enunciated only for those elliptic curves whose $L$-functions have this
modulariy property~; thanks to Wiles and his successors, we now know that they
all have.

Mazur's determination of the possible torsion subgroups (Theorem~18) involves
the study of {\it modular curves}, which are intimately related to modular
forms.

Tunnell's criterion (Theorem~25) is actually an expression for $L(C_\alpha,1)$
in terms of (the ``real period'' of $C_1$ and) the coefficients $c_j(n)$ of
certain modular forms of half-integral weight (cf.~Notation~10).

The role of {\it automorphic forms\/} --- a generalisation of modular forms
--- is even greater in the results of Skinner-Urban (Theorem~29) and in the
proof of the Sato-Tate conjecture (Theorem~33).  It is unlikely to diminish in
the future~: more and more $L$-functions are going to become automorphic,
fulfilling the prophetic vision of Robert Langlands \citer\langlands().

{\eightpoint For a first introduction, apart from Serre's {\it Course}, see
  the book by Koblitz and Knapp's {\it Elliptic curves}.}

\bigbreak
\unvbox\bibbox
\bigbreak
{\obeylines\parskip=0pt\parindent=0pt
Chandan Singh Dalawat
Harish-Chandra Research Institute
Chhatnag Road, Jhunsi
{\pc ALLAHABAD} 211\thinspace019, India
\vskip5pt
\tt dalawat@gmail.com}

\bye